\documentclass[12pt,a4paper]{article}
\usepackage[utf8]{inputenc} 
\usepackage{amsthm,amsmath,amssymb}
\usepackage{graphicx,color}
\usepackage[left=1.5in, right=1in, top=1in, bottom=1in, includefoot, headheight=13.6pt]{geometry}
\usepackage[square, comma, numbers, sort&compress]{natbib}
\usepackage{hyperref}

\usepackage[T1]{fontenc}
\usepackage[toc,page]{appendix}

\usepackage{mathpazo}
\usepackage[hang, small, bf, margin=20pt, tableposition=top]{caption}
\usepackage{booktabs} 
\usepackage{array} 
\usepackage{verbatim} 

\usepackage{epstopdf}

\DeclareMathOperator*{\argmax}{\arg\!\max}

\makeatletter
\newcommand{\mypm}{\mathbin{\mathpalette\@mypm\relax}}
\newcommand{\@mypm}[2]{\ooalign{%
  \raisebox{.1\height}{$#1+$}\cr
  \smash{\raisebox{-.6\height}{$#1-$}}\cr}}
\makeatother

\begin{document}
\centerline{\Large \bf Gaussian asymptotic limits for the $\alpha$-transformation}
\centerline{\Large \bf in the analysis of compositional data}
\vskip 0.2truein
\centerline{\large Yannis Pantazis}
\vskip 0.1truein
\centerline{\it Institute of Applied and Computational Mathematics,} 
\centerline{\it Foundation for Research and Technology - Hellas,} 
\centerline{\it Heraklion 70013, Greece}
\centerline{\it \href{mailto:pantazis@iacm.forth.grm}{pantazis@iacm.forth.gr}}

\vskip 0.2truein
\centerline{\large Michail T. Tsagris}
\vskip 0.1truein
\centerline{\it Department of Economics, University of Crete, } 
\centerline{\it Rethymnon 74100, Greece}
\centerline{\it \href{mailto:mtsagris@uoc.gr}{mtsagris@uoc.gr}}

\vskip 0.2truein
\centerline{and \large Andrew T.A. Wood}
\vskip 0.1truein
\centerline{\it School of Mathematical Sciences, University of Nottingham, }
\centerline{\it Nottingham NG7 2RD, UK}
\centerline{\it \href{mailto:andrew.wood@nottingham.ac.uk}{andrew.wood@nottingham.ac.uk}}

\begin{abstract}
Compositional data consists of vectors of proportions whose components sum to $1$. Such vectors lie in the standard simplex, which is a manifold with boundary.  One issue that has been rather controversial within the field of compositional data analysis is the choice of metric on the simplex. One popular possibility has been to use the metric implied by log-transforming the data, as proposed by Aitchison \cite{ait1983,ait1986}; and another popular approach has been to use the standard Euclidean metric inherited from the ambient space. Tsagris et al. \cite{tsagris2011} proposed a one-parameter family of power transformations, the $\alpha$-transformations, which include both the metric implied by Aitchison's transformation and the Euclidean metric as particular cases. Our underlying philosophy is that, with many datasets, it may make sense to use the data to help us determine a suitable metric.  A related possibility is to apply the  $\alpha$-transformations to a parametric family of distributions, and then estimate $\alpha$ along with the other parameters.  However, as we shall see, when one follows this last approach with the Dirichlet family,  some care is needed in a certain limiting case which arises ($\alpha \rightarrow 0$), as we found out when fitting this model to real and simulated data. Specifically, when the maximum likelihood estimator of $\alpha$ is close to $0$, the other parameters tend to be large.  The main purpose of the paper is to study this limiting case both theoretically and numerically and to provide insight into these numerical findings.
\end{abstract}

\vskip 0.2truein

\vskip 0.2truein

\section{Introduction}

A compositional data vector has non-negative components which sum to a constant, usually taken to be $1$, as we assume here.  Such vectors lie in the standard simplex, $\mathcal{S}\textrm{im}(d)$,  which has the form
\begin{equation}
\mathcal{S}\textrm{im}(d)=\{{\bf x}=(x_1,\ldots , x_D)^\top \in \mathbb{R}^D:\, x_j \geq 0,\,  j=1, \ldots , D;\, \sum_{j=1}^D x_j=1\},
\label{simplex}
\end{equation}
where $d=D-1$.  Compositional datasets appear naturally in many fields including medicine, geology, archeology, biology and economics; see e.g. \cite{ait1986}.  
\medskip

The following two approaches have been widely used in compositional data analysis.  One of these approaches, due to Aitchison \cite{ait1983, ait1986}, is to apply the following transformation, known as the centred log-ratio transformation:
\begin{equation}
w_{j}=\log (x_{j})-D^{-1}\sum_{k=1}^D \log(x_k),  \hskip 0.1truein j=1,\ldots , D.
\label{clr}
\end{equation}
Also, define ${\bf w}=(w_1, \ldots , w_D)^\top$.  Note that, by definition, the components of ${\bf w}$ sum to $0$. Given a sample of compositional vectors ${\bf x}_1, \ldots , {\bf x}_n$, we transform to ${\bf w}_1, \ldots , {\bf w}_n$ using (\ref{clr}). Aitchison then proposed applying standard multivariate methods, such as principal components analysis, to log-transformed compositional vectors, which is tantamount to assuming the standard Euclidean metric is appropriate for  the ${\bf w}_i$. \medskip

Another popular approach is to ``do nothing''; i.e. assume that the Euclidean metric inherited from the ambient Euclidean space $\mathbb{R}^D$ is appropriate for the ${\bf x}_i$, and apply standard multivariate methods directly to the ${\bf x}_i$.  For advocation of this approach see Baxter \cite{baxter1995,baxter2001}, Baxter et al. \cite{baxter2005} and Baxter and Freestone \cite{baxter2006}.\medskip

For some datasets, the log-ratio approach seems preferable, whereas for other datasets, the ``do nothing'' approach seems superior.  Nevertheless, the choice over which method to use, if either, has occasionally become heated and acrimonious.  See the paper by Scealy and Welsh \cite{scealy2014} for a careful and detailed discussion of the issues and history of this debate.  Our view concerning the choice of metric is that {\em it is not possible to come up with a compelling choice for either method based on purely a priori or theoretical grounds}, and that {\em a more pragmatic approach is to make a data-dependent choice of metric}.\medskip

Bearing this in mind, one possibility is to consider a finite-dimensional family of transformations, such as the one-dimensional $\alpha$-transformation family considered in this paper, where each transformation implies a choice of metric. Then use the data to choose the optimal transformation, which in turn determines an optimal metric.  It is not so clear how to do this non-parametrically, but if a parametric family of models, such as the Dirichlet in the present context, is assigned to the transformed compositional data, then we can estimate both the distributional parameters and the transformation parameter simultaneously by maximum likelihood. Using this approach we estimate the parameters of the $\alpha$-transformed compositional data model. In effect, the optimal transformation parameter compels the data to resemble the chosen parametric family as closely as possible.
However, after carrying out this approach with a number of real and simulated datasets we found that in cases where the maximum likelihood estimator of $\alpha$ is close to $0$, there is a tendency for the other (Dirichlet) parameters to be rather large.  See Tables 1-4 in Section 5 for examples of this phenomenon.  This has motivated us to take a close look at the asymptotics as (i) $\alpha \rightarrow 0$ and (ii) the Dirichlet becomes highly concentrated.  This theoretical and numerical investigation is the principal focus of our paper.  \medskip

The phenomenon mentioned in the previous paragraph potentially arises with analogous families of metrics on other manifolds should we wish to select an optimal metric in a data-dependent way.  Consequently, before focusing on compositional data in the remainder of the paper, we briefly discuss the question of choice of metric in the setting of general manifolds in Section 2.\medskip

The structure of this paper is as follows.  In the next section we discuss the choice of metric in the statistical analysis of manifold-valued data from a broader perspective.  Then, in Section 3, we briefly review the $\alpha$-transformation for compositional data.  In Section 4 we present some asymptotics for small $\alpha$ when the $\alpha$-transformation is applied to the Dirichlet family, and in Section 5 we present numerical results, both from simulations and the analysis of real datasets.  These numerical results reflect the theoretical results in the previous section.  Proofs are given in Section 6. \medskip

\section{Statistics on  manifolds and choice of metric}

Since the middle of the last century, statistical analysis of data which naturally lies in a manifold has grown enormously in interest and relevance.  An important early example was R.A. Fisher's paper \cite{fisher1953} on directional statistics with application to paleomagnetic directional data.  Subsequently, the field of directional statistics, in which the sample space consists of the surface of a unit sphere in Euclidean space (or, equivalently, the set of directions or unit vectors) has grown into a well-developed and mature field.  See, for example, the books by Mardia \cite{mardia1972}, Mardia and Jupp \cite{mardia2000} and Fisher et al. \cite{fisher1987}. \medskip

In a related but different direction, there has been great interest in statistical shape analysis from the 1980s onwards.  A popular and mathematically deep approach to statistical shape analysis was initiated by D.G. Kendall; see the monograph Kendall et al. \cite{kendall1999}.  Other books on this field include Small \cite{small1996} and  Dryden and Mardia \cite{dryden1998,dryden2016}.  In the Kendall approach, objects in Euclidean space, usually in two or three dimensions, are represented by the coordinates of a finite set of \textit{labelled landmarks}.  For example, a human face could be represented approximately by landmark coordinates at obvious places including the eyes, ears, nose and mouth.  Then the shape is defined as what remains when translation, scale and orientation have been ``quotiented out''.  Further details may be found in the literature mentioned above.\medskip

Whenever analysing manifold data, one question that always needs to be addressed is what metric to use.  For example, a choice of metric is required if we wish to define a Frech\'et mean on the manifold, which is a generalisation of the least-squares mean in Euclidean space. One important dichotomy is between intrinsic and extrinsic metrics.  In the former, the metric, often a Riemannian metric, is defined ``intrinsically'', without any reference to an ambient space.  In contrast, with extrinsic metrics, the manifold of interest is considered  to be embedded in an ambient space, typically a Euclidean space endowed with the standard Euclidean metric, and the manifold inherits the metric from the ambient space.  From a conceptual point of view, intrinsic metrics are often to be preferred.  However, extrinsic metrics are often easier to work with, both from a practical and theoretical point of view; from a practical point of view because intrinsic distance may be more difficult to calculate and from a theoretical point of view because e.g. asymptotic theory tends to be more difficult when intrinsic metrics are used.  As an example of the difficulties which can arise, see Hotz and Huckemann \cite{hotz2015}, who study the non-standard behaviour which arises when the intrinsic, or arc-length, metric is used for circular data.\medskip

Even when we have made the choice to use an extrinsic metric, the choice of metric may still not be clear. Two examples (which supplement the case of compositional data considered here) in which a one-parameter family of metrics may be worth considering will now be mentioned. \medskip

  In \cite{dryden2009} a family of power metrics on the space of covariances matrices of a given dimension is considered, with a particular focus on diffusion tensor imaging.    In this case, the relevant manifold is the space of positive-definite matrices of given dimension.  There is a close analogy with the $\alpha$-transformation considered here; in their setting the powered quantities are the eigenvalues of the covariance matrices, whereas here the powered quantities are components of compositional vectors.\medskip

The second example relates to statistical shape analysis.  In  \cite{dryden2014} a family of extrinsic metrics is considered for measuring the distance between two shapes described by labelled landmarks.  Here, the  family of metrics consists of powers of eigenvalues of certain symmetric non-negative definite matrices, but in this case the relevant matrices do not have full rank.\medskip

In both of these examples, the question of how to choose the metric arises, as it does in this paper, where the focus is on compositional data.  It is clear that the choice one makes is going to have an effect on the statistical analysis, in some cases a major effect.\medskip

\section{Review of the $\alpha$-transformation}

For the purpose of this paper it is convenient to define the $\alpha$-transformation of a compositional 
vector $\mathbf{x} \in \mathcal{S}\textrm{im}(d)$  by
\begin{eqnarray} \label{stayalpha}
  {\bf u}_\alpha(\mathbf{x})=\left( \frac{x_1^{\alpha}}{\sum_{k=1}^Dx_k^{\alpha}}, \ldots, \frac{x_D^{\alpha}}{\sum_{k=1}^Dx_k^{\alpha}} \right)^T.
\end{eqnarray}
This transformation, which has a slightly different but mathematically equivalent form to that given in \cite{tsagris2011} and \cite{tsagris2016},  is also known as the compositional power transformation \cite{ait1986}.\medskip

Note that the transformation ${\bf u}_\alpha$  given by (\ref{stayalpha}) defines a bijection of the interior of $\mathcal{S}\textrm{im}(d)$
with inverse
\begin{equation}
  {\bf u}_\alpha^{-1}(\mathbf{x})
=\left(\frac{x_1^{1/\alpha}}{\sum_{k=1}^Dx_k^{1/\alpha}},\ldots,\frac{x_D^{1/\alpha}}{\sum_{k=1}^Dx_k^{1/\alpha}}\right).
\label{stayalphaInverse}
\end{equation}
If one excludes the boundary of the simplex, which
corresponds to compositional vectors that have one or more components equal to zero, then the
$\alpha$-transformation (\ref{stayalpha}) and its inverse (\ref{stayalphaInverse}) are well defined
for all $\alpha \in \mathbb{R}$. The motivation for transformation (\ref{stayalpha}) is that the case $\alpha = 0$
corresponds to the log-ratio approach, whereas $\alpha = 1$ corresponds to the "do nothing" approach; see the discussion in the Introduction. We define the case $\alpha = 0$ in
terms of the limit $\alpha \rightarrow 0$; then, rescaling by $\alpha$ as follows, we obtain
\begin{equation*}
\lim_{\alpha \rightarrow 0} \, \alpha^{-1} \left ( D {\bf u}_\alpha({\bf x})-{\bf 1}_D \right ) = \mathbf{w}(\mathbf{x}),
\end{equation*}
where ${\bf 1}_D$ is the $D$-vector of ones and ${\bf w}(\mathbf{x})$ is the centred log-ratio transformation (\ref{clr}). For the case $\alpha = 1$, (\ref{stayalpha}) is just the identity
transformation of the simplex. \medskip

The $\alpha$-transformation (\ref{stayalpha}) leads naturally to a simplicial
distance measure $\Delta_{\alpha}\left({\bf x},{\bf y}\right)$, which we call the
$\alpha$-metric, between observations ${\bf
x},{\bf y} \in \mathcal{S}\textrm{im}(d)$, defined in terms of the Euclidean distance $\| \cdot \|$ between transformed observations, i.e.,
\begin{align*} 
\Delta_{\alpha}\left({\bf x},{\bf y}\right)  = \frac{D}{|\alpha|} \| \mathbf{u}_\alpha(\mathbf{x}) - \mathbf{u}_\alpha(\mathbf{y}) \|
=\frac{D}{\left| \alpha \right|}\left[\sum_{j=1}^D\left(\frac{x_j^\alpha}
{\sum_{k=1}^D x_k^\alpha}-\frac{y_j^\alpha}{\sum_{k=1}^D
y_k^\alpha}\right)^2\right]^{1/2}.
\end{align*}
The special case
\begin{equation*}
\Delta_0\left({\bf x},{\bf y}\right)=
\lim_{\alpha \rightarrow 0} \Delta_\alpha\left({\bf x},{\bf y}\right)
= \left[\sum_{j=1}^D\left(\log{\frac{x_j}{\prod_{k=1}^Dx_k^{1/D}}}-\log{\frac{y_j}{\prod_{k=1}^Dy_k^{1/D}}} \right)^2 \right]^{1/2},
\end{equation*}
is the Euclidean distance on the log-ratio transformed data, whereas
\begin{eqnarray*}
\Delta_1\left({\bf x},{\bf y}\right)=D\left[\sum_{j=1}^D\left(x_j-y_j\right)^2\right]^{1/2}
\end{eqnarray*}
is just Euclidean distance multiplied by $D$.\medskip

The choice of the value of $\alpha$ depends upon the context of the analysis. Maximum likelihood estimation, assuming the transformation (\ref{stayalpha}) or (\ref{stayalphaInverse}) has been applied to a parametric family of distributions such as the Dirichlet, requires the Jacobian determinant of (\ref{stayalpha}), which  is given by (see \cite{tsagris2018} for the derivation)
\begin{eqnarray}  \label{jac}
\left|{\bf J}_{\alpha}\right|=\vert \alpha \vert^d\prod_{j=1}^D\frac{x_j^{\alpha-1}}{\sum_{k=1}^Dx_k^{\alpha}}.
\end{eqnarray}

\section{Gaussian asymptotic limits as $\alpha \rightarrow 0$}

In this section, we examine the $\alpha$-transformation defined in (\ref{stayalpha}), applied to a sample ${\bf x}_1, \ldots , {\bf x}_n$ of compositional data on $\mathcal{S}\textrm{im}(d)$. The $\alpha$-transformation along with the assumption of a $D$-dimensional Dirichlet distribution for the transformed data result in a parametric family which has $D+1$ parameters in total (i.e. one additional parameter). The particular focus in this Section is to see what happens when we apply maximum likelihood estimation to the transformed $D+1$ parameter family in the limit where the Dirichlet becomes highly concentrated and Gaussian and $\alpha \rightarrow 0$.  This question arose when fitting this transformed model to datasets where the MLE of $\alpha$ turns out to be small, while the MLE of the Dirichlet parameters tended to be large.  For numerical evidence of this phenomenon, see Tables 1-4 below.\medskip

\subsection{The Dirichlet case}


Recall that the Dirichlet density on the simplex at a point ${\bf v}=(v_1, \ldots , v_D)^\top \in \mathcal{S} \textrm{im}(d)$ is 
\[
f_{\pmb \gamma}({\bf v}) =\frac{\Gamma (\gamma_+)}{\prod_{j=1}^D \Gamma(\gamma_j)}\prod_{j=1}^D v_j^{\gamma_j-1},
\]
where $\gamma_+=\sum_{j=1}^D \gamma_j$. Suppose that
\[
{\bf u}_i=(u_{i1}, \ldots , u_{iD})^\top \sim \textrm{Dirichlet}(b_1,\ldots ,b_D), \hskip 0.2truein i=1,\ldots , n,
\]
where
\[
u_{ij}=\frac{x_{ij}^\alpha}{\sum_{k=1}^D x_{ik}^\alpha}, \hskip0.1truein i=1, \ldots, n, \hskip 0.1truein j=1,\ldots , D.
\]
Then the log-likelihood which results when the $\alpha$-transformation is applied to ${\bf x}_1,...,{\bf x}_n$ is given by
\[
\bar{l}(\alpha, b_1, \ldots , b_D) = \sum_{i=1}^n \big[ \log f_{\bf b} ({\bf u}_i) + \log \left|{\bf J}_{\alpha,i}\right| \big]
\]
which can be rewritten after expansion and re-organization as
\begin{equation} \label{lik}
\begin{aligned}
\bar{l}(\alpha, b_1, \ldots , b_D) &= n \log(\Gamma(b_+)) - n\sum_{j=1}^D \log(\Gamma(b_j)) + n d\log(|\alpha|) \\
&\hskip 0.3truein + \sum_{i=1}^n\sum_{j=1}^D (\alpha b_j-1) \log(x_{ij}) 
- b_+ \sum_{i=1}^n   \log\left(\sum_{j=1}^D x_{ij}^\alpha \right),
\end{aligned}
\end{equation}
where $b_+=\sum_{j=1}^D b_j$.
Also, define 
\begin{equation}
y_{ij}=\log (x_{ij}), \hskip 0.2truein i=1, \ldots , n; \hskip 0.1truein j=1, \ldots , D.
\label{ydef}
\end{equation}

\vskip 0.2truein

We now state our first result about the asymptotic behavior of the log-likelihood.
\vskip 0.2truein

\noindent {\bf Proposition 1}.  \textit{Define} $y_{ij}$ \textit{as in (\ref{ydef}) and}
\begin{equation}
b_j=\frac{b}{\alpha^2}(1+\alpha c_j), \hskip 0.2truein j=1,\ldots , D,
\label{b_j}
\end{equation}
\textit{where}
\begin{equation}
\sum_{j=1}^D c_j=0,
\label{sum_c}
\end{equation}
\textit{and define}
\[
\ell(\alpha, b, c_1, \ldots , c_D)=\bar{\ell}\left(\alpha, \frac{b}{\alpha^2}(1+\alpha c_1), \ldots , \frac{b}{\alpha^2}(1+\alpha c_D)\right).
\]
\textit{Then as $\alpha \rightarrow 0$,  the log-likelihood $\ell(\alpha, b,c_1,\ldots, c_D)$ for $b$ and the $c_j$ satisfies}
\begin{align}
\ell(\alpha, b,c_1, \ldots , c_D) &=\frac{nd}{2}\log\left (\frac{b}{2\pi} \right) -\frac{b}{2} \sum_{i=1}^n \sum_{j=1}^D \left ( y_{ij} - \bar{y}_{i+} -c_j \right)^2\nonumber \\
 &\hskip 1.0truein +  n\mathcal{C}_{0}+\alpha \sum_{i=1}^n \mathcal{C}_{1,i} + O(\alpha^2),
\label{log-like-limit}
\end{align}
\textit{where} $\bar{y}_{i+}=D^{-1}\sum_{j=1}^D y_{ij}$,\ \  $\bar{y}_{++}=(nD)^{-1}\sum_{i=1}^n \sum_{j=1}^D y_{ij}$,
\begin{equation}
\mathcal{C}_0=-\frac{1}{2} \log(D) -D\bar{y}_{++}, \hskip 0.2truein \mathcal{C}_{1,i}=-\frac{b}{6}\left ( D \hat{\kappa}_{3,i} +\sum_{j=1}^D c_j^3 \right ),
\label{mathcalC}
\end{equation}
\textit{and}
\begin{equation}
\hat{\kappa}_{1,i}\equiv \bar{y}_{1,i}=D^{-1}\sum_{j=1}^D y_{ij}, \hskip 0.2truein \hat{\kappa}_{2,i} = \left (D^{-1} \sum_{j=1}^D y_{ij}^2\right)-\hat{\kappa}_{1,i}^2
\label{kappa12}
\end{equation}
\textit{and}
\begin{equation}
\hat{\kappa}_{3,i}=\left (D^{-1}\sum_{j=1}^D y_{ij}^3\right) -3 \hat{\kappa}_{1,i}\left (D^{-1} \sum_{j=1}^D y_{ij}^2\right) + 2 \hat{\kappa}_{1,i}^3
\label{kappa3}
\end{equation}
\textit{are the first three sample cumulants of the} $y_{i1}, \ldots , y_{iD}$.

\vskip 0.2truein

\noindent Proposition 1 is proved in Section 6. We present several remarks below concerning Proposition 1.
\vskip 0.2truein

\noindent \textit{Remark 1}.  When $\alpha$ is small, the MLE obtained by maximising (\ref{log-like-limit}) is given by
\begin{equation} \label{c:hat:eq}
\hat{c}_j =n^{-1}\sum_{i=1}^n (y_{ij}- \bar{y}_{i+})+ O(\alpha)=\bar{y}_{+j}-\bar{y}_{++}+O(\alpha), \hskip 0.2truein 1 \leq j \leq D,
\end{equation}
where $\bar{y}_{+j}=n^{-1}\sum_{i=1}^n y_{ij}$,
and
\begin{equation} \label{b:hat:eq}
\hat{b}=\left \{\frac{1}{nd} \sum_{i=1}^n \sum_{j=1}^D \left (y_{ij}-\bar{y}_{+j}-\bar{y}_{i+}+\bar{y}_{++}\right )^2\right \}^{-1} +O(\alpha)
\end{equation}

\vskip 0.2truein

\noindent \textit{Remark 2}.  The limit distribution is logistic normal with isotropic covariance structure (see Remark 3 below for clarification).  The simple covariance structure in the limit as $\alpha \rightarrow 0$ is a  consequence of the restricted nature of the covariance structure of the Dirichlet distribution.  Note that the limit distribution has the same number of parameters as the $D$-dimensional Dirchlet distribution, namely $D$.

\vskip 0.2truein

\noindent \textit{Remark 3}.  Although the likelihood in (\ref{log-like-limit}) clearly has a Gaussian character, it may not be immediately clear how it arises, given that when we sum
$y_{ij}-\bar{y}_{i+}$ over $j$ we get $0$.  To clarify, suppose we have independent Gaussian variables $Z_j \sim N(\mu_j,\sigma^2)$, $j=1, \ldots , D$, where it is assumed that $\sum_{j=1}^D \mu_j=0$.  Define $\bar{Z}=D^{-1}\sum_{j=1}^D Z_j$ to be the sample mean.  The covariance matrix of $(Z_1-\bar{Z}, \ldots , Z_D-\bar{Z})^\top$ is ${\pmb \Sigma}=\sigma^2 \left ({\bf I}_D-D^{-1}{\bf 1}_{D}{\bf 1}_{D}^\top \right)$, where ${\bf I}_D$ is the $D \times D$ identity matrix and ${\bf 1}_D$ is the $D$-vector of ones.  The Moore-Penrose inverse of ${\pmb \Sigma}$ is given by
\[
{\pmb \Sigma}^{-}=\frac{1}{\sigma^2}\left ( {\bf I}_D-D^{-1} {\bf 1}_D{\bf 1}_D^\top  \right )
\]
and the relevant quadratic form is given by
\begin{align}
&(Z_1-\bar{Z}-\mu_1, \ldots , Z_D-\bar{Z}-\mu_D){\pmb \Sigma}^{-}(Z_1-\bar{Z}-\mu_1, \ldots , Z_D-\bar{Z}-\mu_D)^\top \nonumber\\
&=\frac{1}{\sigma^2}\sum_{j=1}^D (Z_j-\bar{Z}-\mu_j)^2-\frac{1}{D\sigma^2}\left (\sum_{j=1}^D  (Z_j-\bar{Z}-\mu_j)  \right )^2 \nonumber \\
&=\frac{1}{\sigma^2}\sum_{j=1}^D (Z_j-\bar{Z}-\mu_j)^2,
\label{crime_and_punishment}
\end{align}
because we have assumed that $\mu_1+\cdots +\mu_D=0$.  Note that if we put $b=1/\sigma^2$, $c_j=\mu_j$ and $Z_j=y_{ij}$, $1 \leq j \leq D$, then
(\ref{crime_and_punishment}) is equal to
\[
b\sum_{j=1}^D \left ( y_{ij}-\bar{y}_{i+}-c_j  \right )^2,
\]
which agrees with the inner sum of the second term in (\ref{log-like-limit}).

\vskip 0.2truein

\noindent \textit{Remark 4}.  The assumption (\ref{sum_c}) only plays a cosmetic role in the proof provided that we are prepared to allow $b$ and the $c_j$ to depend weakly on $\alpha$.  Specifically, suppose that $\bar{c}_+=D^{-1} \sum_{j=1}^D c_j \neq 0$.  Then define
\[
b^\ast=b(1+\alpha \bar{c}_+) \hskip 0.2truein \hbox{and} \hskip 0.2truein c_j^\ast =\frac{c_j-\bar{c}_+}{1+\alpha \bar{c}_+}, \hskip 0.2truein j=1, \ldots , D.
\]
It is easily checked that
\[
b_j\equiv \frac{b}{\alpha^2}(1+\alpha c_j)=\frac{b^\ast}{\alpha^2}(1+\alpha c_j^\ast) \hskip 0.2truein \hbox{and} \hskip 0.2truein \sum_{j=1}^D c_j^\ast = 0,
\]
 and also $b^\ast$ and the $c_j^\ast$ only depend weakly on $\alpha$ in the sense that $b^\ast \rightarrow b$ and $c_j^\ast \rightarrow c_j-\bar{c}_+$, $1 \leq j \leq D$, as $\alpha \rightarrow 0$.

\vskip 0.2truein

\noindent \textit{Remark 5}. Note that the $b_j$ defined in (\ref{b_j}) have a rather special structure in that they are asymptotically equal as $\alpha \rightarrow 0$.  A natural question is: what happens to the asymptotics in Proposition 1 if the $b_j$ do not coalesce as $\alpha \rightarrow 0$, specifically if the $b_j$ are of the form $b_j=\gamma_j/\alpha^2$.  It turns out that, in contrast to Proposition 1, when the $b_j$ do not coalesce the log-likelihood does not converge to a finite limit as $\alpha \rightarrow 0$, as we shall see later.

\subsection{The role of asymptotic normality and the general case}

The asymptotic regime considered in the previous section is quite limited because the $b_j$ in (\ref{b_j}) satisfy
\[
\frac{b_j}{\sum_{k=1}^D b_k} \rightarrow \frac{1}{D} \hskip 0.2truein \textrm{as} \hskip 0.1truein \alpha \rightarrow 0.
\]
However, it seems that we need this strong assumption to hold if we are to have convergence to a finite log-likelihood as indicated in Proposition 1.  One important point we have not discussed yet is the fact that the Dirichlet distribution in Proposition 1 is converging to a multivariate normal.  Consider the following elementary result.\medskip

\noindent \textbf{Proposition 2}.  \textit{Suppose} ${\bf u}^{(\alpha)} \sim \textrm{Dirchlet}({\pmb \gamma}/\alpha^2)$, where ${\pmb \gamma} =(\gamma_1, \ldots , \gamma_D)^\top$.  \textit{Then as} $\alpha \rightarrow 0$,
\[
\frac{1}{\alpha}({\bf u}^{(\alpha)} - \bar{\pmb \gamma}) \rightarrow^{\hskip -0.13truein d} \mathcal{N}_D({\bf 0}_p, {\pmb \Sigma}),
\]
\textit{where} ${\pmb \Sigma}=\gamma_+^{-1}\left \{\textrm{diag}(\bar{\pmb \gamma})-\bar{\pmb \gamma}\bar{\pmb \gamma}^\top\right \}$,
$\bar{\pmb \gamma}=\gamma_{+}^{-1}{\pmb \gamma}$ \textit{and} $\gamma_+=\sum_{j=1}^D \gamma_j$.  \textit{Moreover, the convergence of densities also occurs uniformly in growing balls of the form } ${\bf v} \in \mathcal{B}_{\alpha^{-1+\epsilon}}(0)$ \textit{as} $\alpha \rightarrow 0$ \textit{for any fixed} $\epsilon \in (0,1)$, \textit{in the sense that, on the set} $\sum_{j=1}^D v_j =0$,
\[
\alpha^d f_{{\pmb \gamma}/\alpha^2}(\bar{\pmb \gamma}+ \alpha{\bf v})\rightarrow \frac{1}{(2\pi)^{d/2}}\left (\frac{\gamma_+^{d}}{\prod_{j=1}^D \bar{\gamma}_j}\right )^{1/2} \exp \left (-\frac{1}{2}\sum_{j=1}^D \frac{\gamma_+ v_j^2}{ \bar{\gamma}_j}\right )
\]
\textit{uniformly for} ${\bf v} \in \mathcal{B}_{\alpha^{-1+\epsilon}}(0)$, \textit{where in the line above,} $f$ \textit{is the Dirichlet density.}
\medskip

The proof of Proposition 2 is similar to that of Proposition 1 but is somewhat simpler and so  is omitted.
\medskip

\noindent We now consider the asymptotic behaviour of $\bf y$ as $\alpha \rightarrow 0$ in the context of a Dirichlet distribution with parameters which do not coalesce as in Proposition 1.\medskip

\noindent \textbf{Corollary 1}.  \textit{Under the assumptions of Proposition 2,}
\begin{equation}
{\bf y} - \bar{y}{\bf 1}_D -\frac{D}{\alpha}\left (\bar{\pmb \gamma}-D^{-1} {\bf 1}_D \right) \rightarrow^{\hskip -0.13truein d} \mathcal{N}_D({\bf 0}_p, D^2 {\pmb \Sigma}),
\label{threescoreandten}
\end{equation}
\textit{where} ${\bf y}=(y_1, \ldots, y_D)^\top$, $y_j=\log(x_j)$ and $\bar{y}=D^{-1}\sum_{j=1}^D y_j$,
$u_j=x_j^{\alpha}/\sum_{k=1}^D x_k^\alpha$ \textit{and} ${\bf u} =(u_1, \ldots , u_D)^\top \sim  \text{Dirichlet}(\gamma_1/\alpha^2, \ldots , \gamma_D/\alpha^2)$.

\vskip 0.2truein

\noindent \textbf{Proof of Corollary 1}.  If $u_j=x_j^\alpha/\sum_{k=1}^D x_k^\alpha$, then elementary calculations show that as $\alpha \rightarrow 0$,
\begin{equation}
u_j=\frac{1}{D} +\frac{\alpha}{D}(y_j-\bar{y}) +O(\alpha^2).
\label{uj2}
\end{equation}
where $y_j=\log (x_j)$ and $\bar{y}=D^{-1}\sum_{j=1}^D y_j$.
 Consequently, from (\ref{uj2}),
\begin{equation}
{\bf y}-\bar{y}{\bf 1}_D =\frac{D}{\alpha}({\bf u}-D^{-1}{\bf 1}_D)+O(\alpha),
\label{bit2}
\end{equation}
and so Corollary 1 follows because, from Proposition 2, 
${\bf u} \approx \mathcal{N}_D(\bar{\pmb \gamma}, \alpha^2 {\pmb \Sigma})$.  \hfill $\square$
\vskip 0.2truein
So for ${\bf u}$ to attain a general limiting mean vector $\bar{\pmb \gamma}$ in the limit as $\alpha \rightarrow 0$, it is necessary for the expectation vector of ${\bf y}-\bar{y}{\bf 1}_D$ to go to infinity as $\alpha \rightarrow 0$.\medskip

Finally, we go to a more general case in which we do not assume that ${\bf u}\equiv {\bf u}^{(\alpha)}$ is Dirichlet but just assume that as $\alpha \rightarrow 0$,
\[
\alpha^{-1} ({\bf u}-{\pmb \eta}) \rightarrow^{\hskip -0.13truein d} \mathcal{N}_D({\bf 0}_D,  {\pmb \Omega}),
\]
Since the ${\bf u}$ are defined on the simplex $\mathcal{S}\textrm{im}(d)$ the vector ${\pmb \eta}=(\eta_1, \ldots , \eta_D)^\top$ has non-negative components which satisfy $\sum_{j=1}^D \eta_j=1$ and the covariance matrix ${\pmb \Omega}$ has an eigenvector ${\bf 1}_D$ with corresponding eigenvalue $0$.  In this case, we obtain the same limiting result as in (\ref{threescoreandten}), but with the more general covariance matrix ${\pmb \Omega}$ replacing ${\pmb \Sigma}$, and with the same requirement that the mean vector of ${\bf y}-\bar{y}{\bf 1}_D$ goes to infinity if we wish for ${\bf u}$ to obtain a general limiting mean vector.

\section{Numerical results}
We first show the asymptotic results of the previous section on synthetic data and then apply the $\alpha$-transformation on real datasets.

\subsection{Numerical Simulations}
\begin{figure}[htb]
\begin{center}
\includegraphics[width=\textwidth]{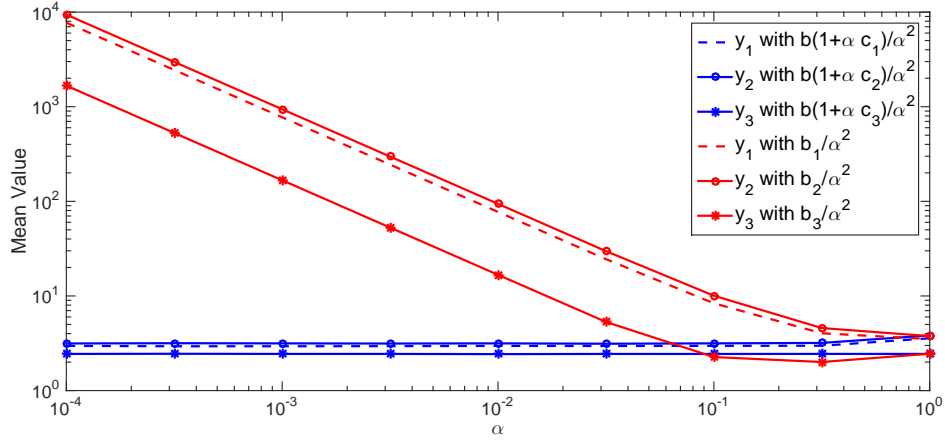}
\caption{Asymptotic comparison as $\alpha\to 0$ between $Dirichlet\left(\frac{b}{\alpha^2}({\bf 1}_D+\alpha {\bf c})\right)$ (blue lines) and $Dirichlet\left(\frac{\bf b}{\alpha^2}\right)$ (red lines).}
\label{num:sim:asym}
\end{center}
\end{figure}

As discussed earlier, a finite log-likelihood is obtained under the assumption in (\ref{b_j}) for $\bf{b}$. In contrast, when the Dirichlet's parameters do not converge to the same value as $\alpha\to 0$ the log-likelihood is expected to be infinite. Figure~\ref{num:sim:asym} presents this behaviour as it is encoded by the asymptotic limit for the mean value ${\bf y}-\bar{y}$. Blue lines correspond to the former case with $b=1$ and ${\bf c}=[0.1,0.3,-0.4]^T$. As expected from Proposition 1, the mean value of ${\bf y}-\bar{y}$ converges to a finite value which is however different at each coordinate due to the different values of ${\bf c}$'s elements. Red lines correspond to the case with ${\bf b}=[1.1,1.3,0.6]^T$ where the mean value diverges as implied by Corollary 1. Evidently, the rate of divergence is inverse proportional to $\alpha$.

\vskip 0.2truein

\noindent \textit{Remark 6}. The estimation of ${\bf y}$ requires the exponentiation of ${\bf u}$ with $1/\alpha$ which leads to numerical instability in the general case (i.e., the red lines) when $\alpha$ is below $0.001$. In order to ensure numerical stability, we calculate ${\bf y}$ with the following equivalent formula 
\[
y_j = \frac{1}{\alpha}\log\frac{u_j}{u_{j^*}} - \log\left(1+\sum_{k=1,k\neq j^*}^D \left(\frac{u_k}{u_{j^*}}\right)^{\frac{1}{\alpha}} \right)
\]
where $j^*=\argmax_j b_j$.

\subsection{Real data examples}
Three real datasets are used for illustration purposes. 
\begin{itemize}
\item \textit{Mammals}. This dataset is taken from Hartigan \cite{hartigan1975} and describe the percentage composition of 24 mammals’ milk on the basis of 5 different constituents (water, protein, fat, lactose and ash). A 4-group solution has been considered as the optimal grouping of such data.

\item \textit{East Bay Clams}. From the many colonies of clams in East Bay, 20 were selected at random and from each a sample of clams was taken \cite{ait1986}. Each sample was sieved into three size ranges, large, medium, small; then each size range was sorted by shell colour, dark, light. For each colony the proportions of clams in each colour-size combination was estimated and the corresponding compositions are recorded. 

\item \textit{OECD}. The dataset contains percentages of the labour force and the per capita income of 20 European OECD countries in 1960. The data can be downloaded from the \href{http://lib.stat.cmu.edu/DASL/Datafiles/oecdat.html}{DASL} library.

\item \textit{GRTA}. This dataset contains information on the number of Greek road traffic accidents and persons injured from January 2010 up to August 2017 on a monthly basis. Three quantities of interest are measured, Killed, seriously injured and slightly injured.                                                       
\end{itemize}

\begin{figure}[t]
\centering
\begin{tabular}{cc}
\includegraphics[scale = 0.3]{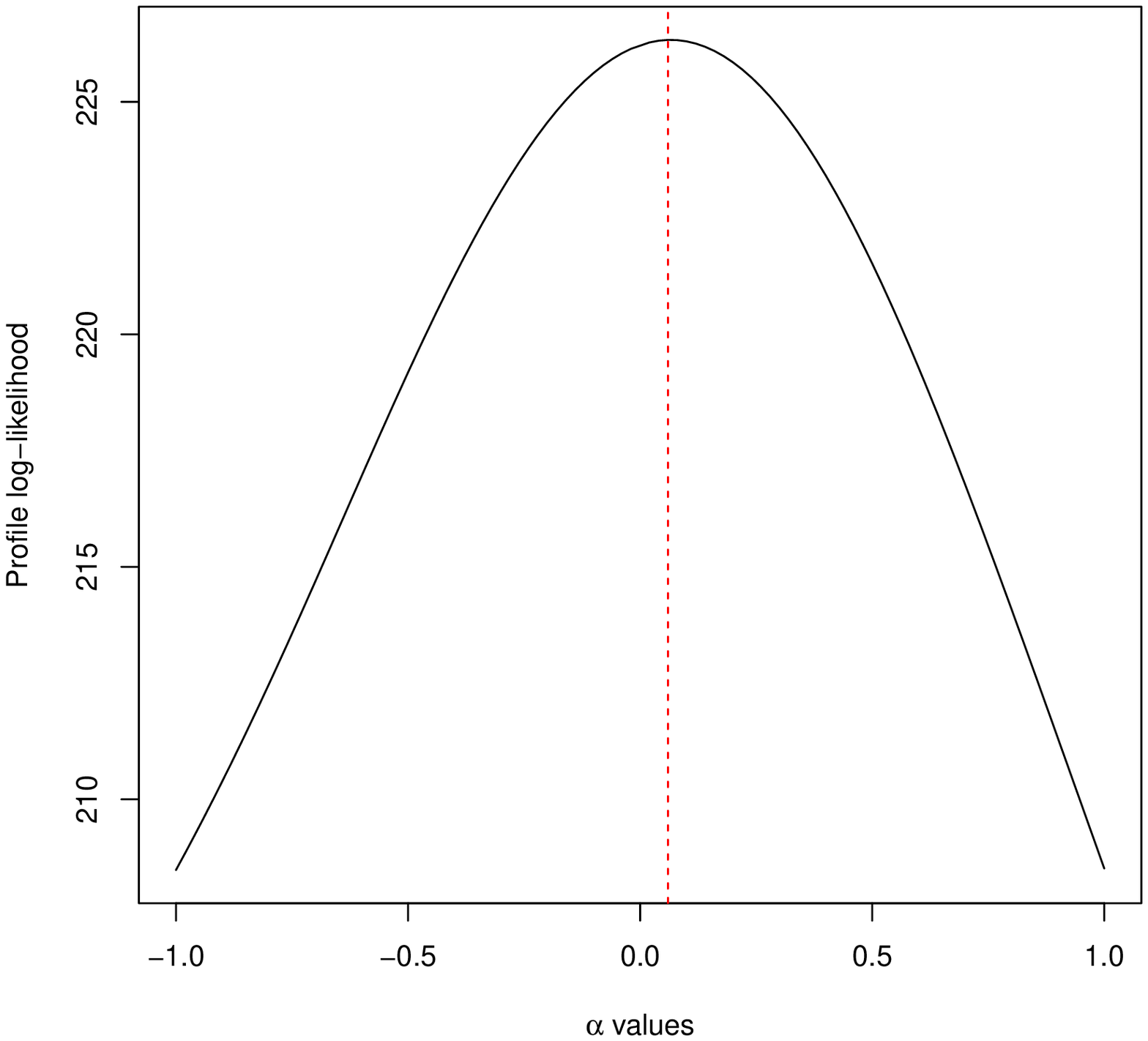} &
\includegraphics[scale = 0.3]{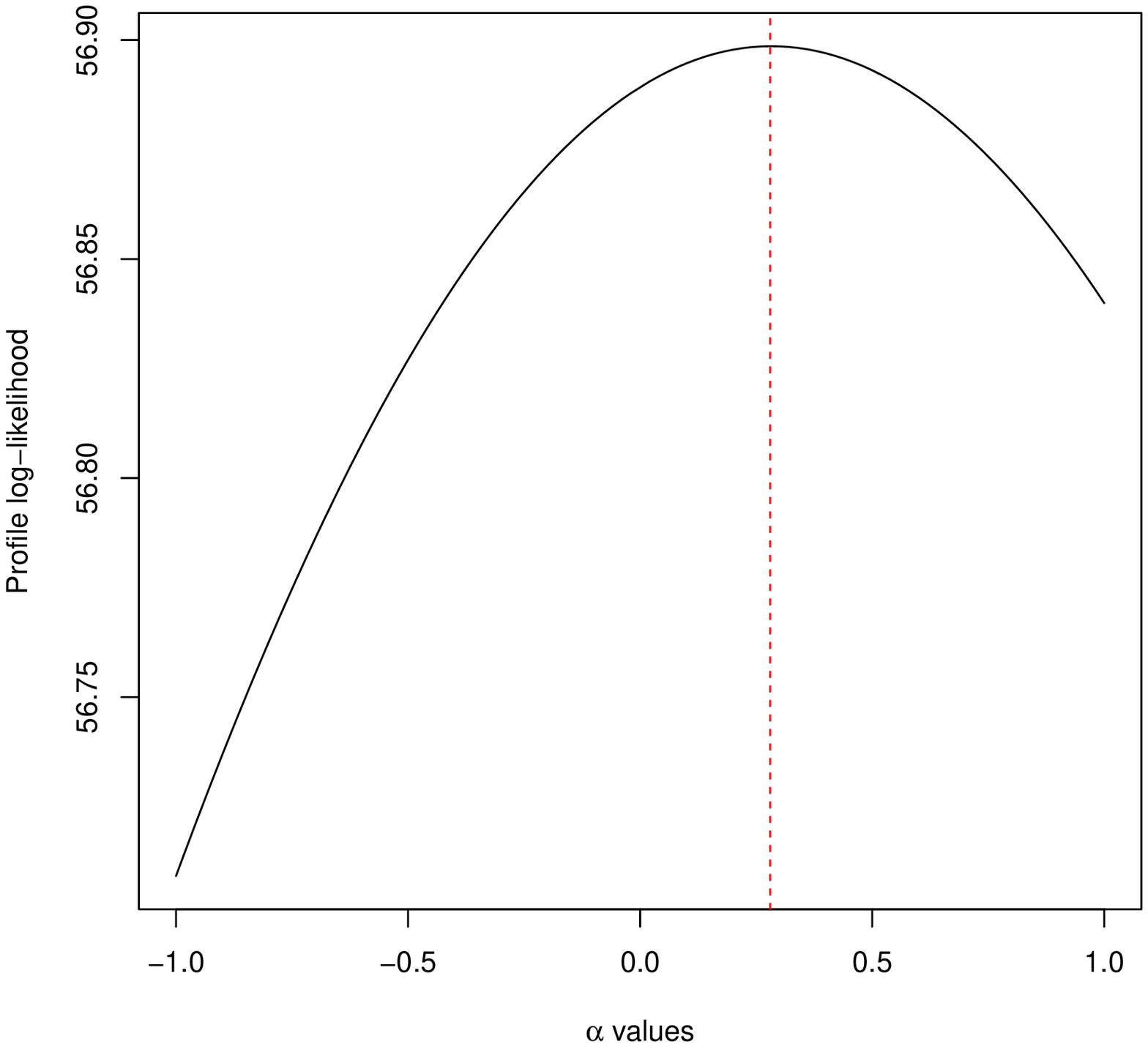}   \\
(a)  &  (b)   \\
\includegraphics[scale = 0.3]{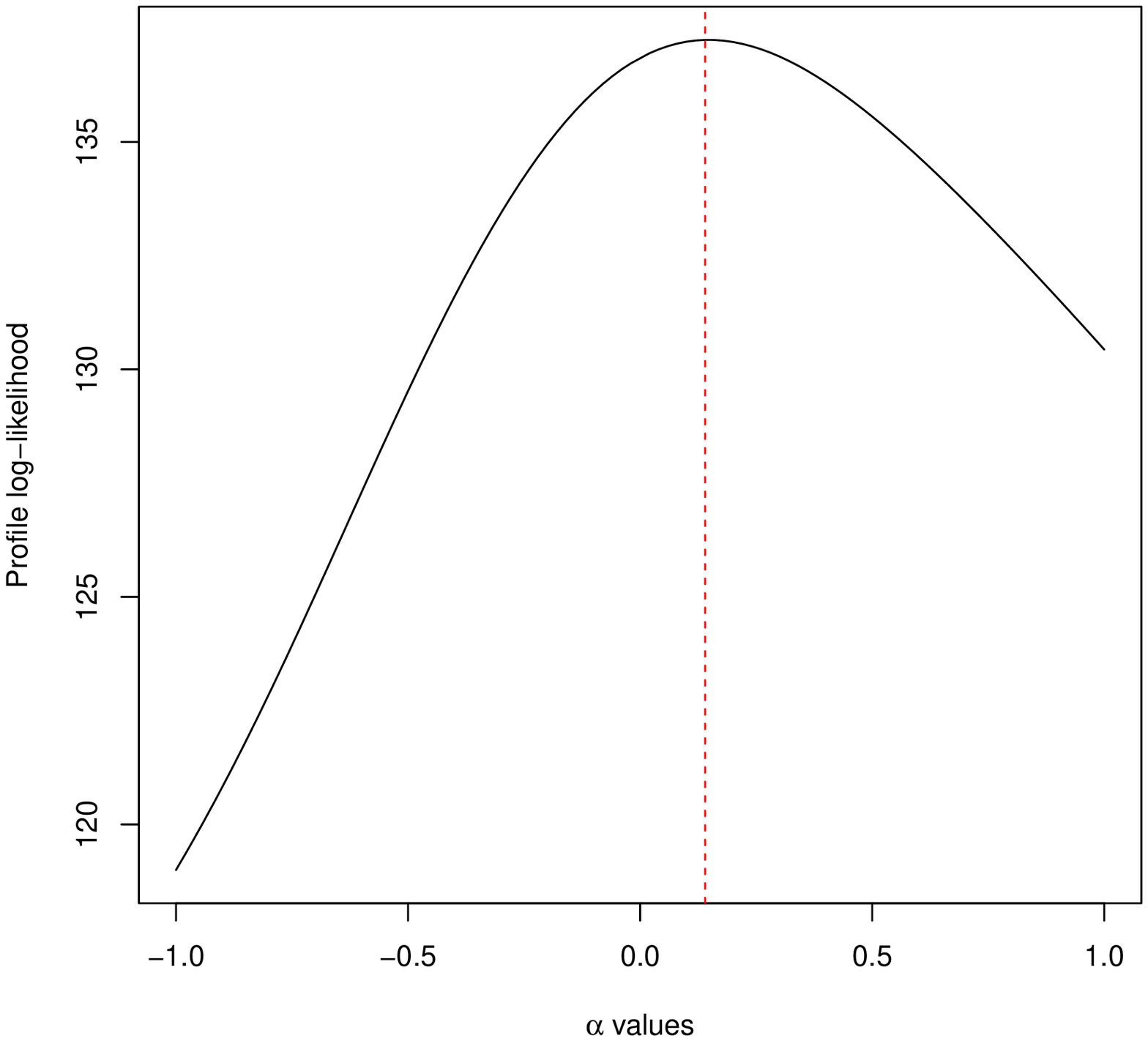}  &
\includegraphics[scale = 0.3]{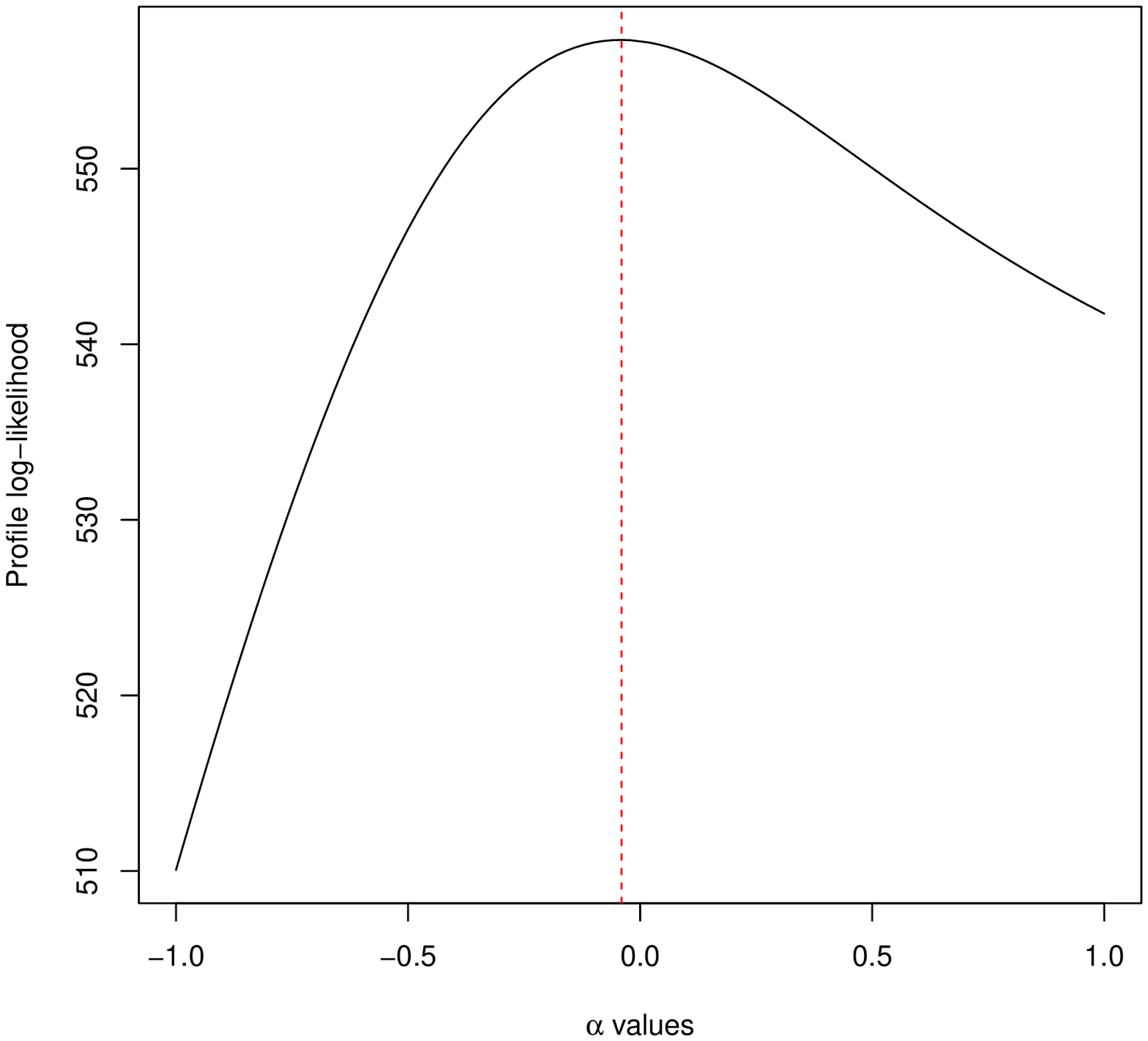}  \\
(a)  &  (b)   \\
\end{tabular}
\caption{Profile log-likelihood as a function of $\alpha$ for the (a) Mammals, (b) East Bay Clams, (c) OECD and (d) GRTA datasets.}
\label{alpha_plots}
\end{figure}

For each of these three datasets, we apply the $\alpha$-transformation (\ref{stayalpha}) and maximise the  Dirichlet log-likelihood. This is the same as maximizing the profile log-likelihood of $\alpha$ presented in (\ref{lik}). For a range of values of $\alpha$ we apply the $\alpha$-transformation (\ref{stayalpha}) to the compositional data and maximize the sum of the Dirichlet log-likelihood in (\ref{lik}) and the logarithm of the Jacobian in (\ref{jac}) of the $\alpha$-transformation given in (\ref{stayalpha}). Figure \ref{alpha_plots} presents these profile log-likelihoods. Tables \ref{mammals}--\ref{grta} show the resulting estimates of the Dirichlet parameters obtained from the MLE applied to the $\alpha$-transformed data and based on the two asymptotic forms. Asymptotic 1 corresponds to $\textrm{Dirichlet}\left(b({\bf 1}_D+\alpha {\bf c})/\alpha^2\right)$ model with $b$ and ${\bf c}$ being estimated from (\ref{b:hat:eq}) and (\ref{c:hat:eq}), respectively while Asymptotic 2 corresponds to the general case $\textrm{Dirichlet}\left({\bf b}/\alpha^2\right)$. The estimation of the parameter vector ${\bf b}$ is based again on MLE. Therefore the results for Asymptotic 2 are expected to be similar with the results when direct MLE is applied. 

\begin{table}[!ht]
\caption{Estimates of the parameters for the Mammals dataset. The columns of the table give estimates of  parameters $b_j$ defined in (\ref{b_j}), but labelled by the component name rather than $j$. The optimal value for $\alpha$ is $0.06$. \label{mammals}}
\centering
\begin{tabular}{lccccc}  
\noalign{\smallskip}\hline\noalign{\smallskip}
Method & water    &   protein  &  fat      &  lactose  &  ash  \\ \hline
Direct MLE &  793.908  &   671.054  &  680.629  &  663.055  &  604.305 \\  
Asymptotic 1 & 782.297 & 668.362 &  677.964  &  660.235 & 597.352 \\
Asymptotic 2 &  793.673 & 670.856 & 680.428  &  662.859  &  604.127  \\  \hline
\end{tabular}
\end{table}

\begin{table}[!ht]
\caption{Estimates for the East Bay Clams dataset. The columns of the table give estimates of  parameters $b_j$ defined in (\ref{b_j}), but labelled by component name rather than $j$. The optimal value for $\alpha$ is $0.28$. \label{east}}
\centering
\begin{tabular}{lccc}  
\noalign{\smallskip}\hline\noalign{\smallskip}
Method  &  dl       &  dm       &  ds   \\ \hline
Direct MLE     &  232.218  &  224.377  &  222.197 \\ 
Asymptotic 1 & 231.568  & 223.800 & 221.592 \\
Asymptotic 2 & 231.990 & 224.156 & 221.979 \\  \hline
\end{tabular}
\end{table}

\begin{table}[t]
\caption{Estimates for the OECD dataset. The columns of the table give estimates of the parameters $b_j$ defined in (\ref{b_j}), but labelled by component name rather than $j$. The optimal value for $\alpha$ is $0.14$. \label{oecd}}
\centering
\begin{tabular}{lcccc}  
\noalign{\smallskip}\hline\noalign{\smallskip}
Method & PCINC      &  AGR     &  IND   &  SER  \\  \hline
Direct MLE    &  212.965  &  129.782 & 143.360  &  136.907  \\  
Asymptotic 1 &  199.034 & 124.925 & 139.824 &  132.928 \\
Asymptotic 2 & 212.669  & 129.601 & 143.161 & 136.716 \\  \hline
\end{tabular}
\end{table}

\begin{table}[t]
\caption{Estimates for the GRTA dataset. The columns of the table give estimates of the parameters $b_j$ defined in (\ref{b_j}), but labelled by component name rather than $j$. The optimal value for $\alpha$ is $-0.04$. \label{grta}}
\centering
\begin{tabular}{lccc}  
\noalign{\smallskip}\hline\noalign{\smallskip}
Method        &  Killed    &  Seriously injured   &  Slightly injured  \\  \hline
Direct MLE    &  28366.65  &  28109.32  &  25420.99  \\  
Asymptotic 1  &  28305.94  &  28057.80  &  25320.62 \\
Asymptotic 2  &  28366.36  &  28109.03  &  25420.73 \\  \hline
\end{tabular}
\end{table}

\vskip 0.2truein

\section{Proofs} 

\noindent \textbf{Proof of Proposition 1}.  Let us focus on the log-likelihood for a single observation $i$.  At the end of the proof we will sum over $i$ to obtain the log-likelihood (\ref{log-like-limit}).\medskip

\noindent After re-organization, the log-likelihood for observation $i$  may be written
\begin{align}
\bar{\ell}_i(\alpha,b,,c_1, \ldots , c_D)&=\log \left \{\Gamma \left (\frac{bD}{\alpha^2} \right )\right \}
-\sum_{j=1}^D \log \left \{\Gamma\left (\frac{b}{\alpha^2}(1+\alpha c_j)\right)\right \}
\nonumber \\
& \hskip 0.5truein +d\log (|\alpha|)
 +\sum_{j=1}^D \left \{  \frac{b}{\alpha}(1+\alpha c_j)-1\right \} \log(x_{ij})\nonumber \\
& \hskip 1.0truein -\sum_{j=1}^D \left \{ \frac{b}{\alpha^2}(1+\alpha c_j) \right \}\log \left (\sum_{k=1}^D x_{ik}^\alpha \right ).
\label{log-likelihood}
\end{align}
We now make use of Lemma 1 and Lemma 2 which are stated and proved below.
Noting that
\[
\sum_{j=1}^D \left \{  \frac{b}{\alpha}(1+\alpha c_j)-1\right \} \log(x_{ij})=\frac{bD}{\alpha}\hat{\kappa}_{1,i}
+b\sum_{j=1}^D c_j y_{ij}-D \hat{\kappa}_{1,i},
\]
and substituting (\ref{lemma1}) and (\ref{lemma2}) into (\ref{log-likelihood}), we obtain, as $\alpha \rightarrow 0$,
\begin{align}
\ell_i(\alpha, b,c_1, \ldots , c_D)&=\bar{\ell}_i\left(\alpha, \frac{b}{\alpha^2}(1+\alpha c_1), \ldots , \frac{b}{\alpha^2}(1+\alpha c_D)\right)\nonumber \\
&=\frac{d}{2}\log (b) -\frac{d}{2}\log(2\pi) - \frac{b}{2}\sum_{j=1}^D \left (y_{ij} -\bar{y}_{i+}-c_j \right )^2 \nonumber \\ 
& \hskip 1.0truein -\frac{1}{2}\log(D) -\frac{\alpha b}{6}\left ( D \hat{\kappa}_{3,i} +\sum_{k=1}^D c_k^3 \right ) +O(\alpha^2)\nonumber \\
&=\frac{d}{2}\log \left (\frac{b}{2\pi}\right)-\frac{b}{2}\sum_{j=1}^D \left (y_{ij} -\bar{y}_{i+}-c_j \right )^2 \nonumber \\
& \hskip 1.0truein + \mathcal{C}_0 +\alpha \mathcal{C}_{1,i}+O(\alpha^2),
\label{likei}
\end{align}
where $\mathcal{C}_0$ and $\mathcal{C}_{1,i}$ are defined in (\ref{mathcalC}). 
Finally, sum (\ref{likei}) over $i$ to obtain (\ref{log-like-limit}).  \hfill$\square$

\vskip 0.2truein
\noindent Lemma 1 and Lemma 2 are now stated and proved.
\vskip 0.2truein

\noindent {\bf Lemma 1}.  \textit{With the $b_j$ defined as in (\ref{b_j}) with the $c_j$ subject to (\ref{sum_c}), the following result holds: as $\alpha \rightarrow 0$,}
\begin{align}
\log \left \{\Gamma \left (\sum_{j=1}^D b_j \right )\right \} &- \sum_{j=1}^D \log \left \{\Gamma(b_j)\right \} \nonumber \\
& =\frac{b D\log(D)}{\alpha^2}-d\log(|\alpha|) +\frac{d}{2}\log(b)  -\frac{1}{2}\log(D)
-\frac{d}{2}\log(2\pi) \nonumber \\
& \hskip 1.0truein -\frac{b}{2}\sum_{j=1}^D c_j^2 - \frac{\alpha b}{6}\sum_{j=1}^D c_j^3+O(\alpha^2)
\label{lemma1}
\end{align}

\vskip 0.2truein

\noindent \textbf{Proof of Lemma 1}.  This involves no more than repeated application of Stirling's formula
\[
\log \{\Gamma(x)\}=\left (x-\frac{1}{2}\right)\log(x) - x+\frac{1}{2}\log(2\pi) +O(x^{-1}), \hskip 0.2truein x \rightarrow \infty.
\]
Making use of (\ref{sum_c}) and applying Stirling's formula gives
\begin{align}
\log \left \{ \Gamma \left ( \sum_{j=1}^D b_j\right) \right \}&=\log \left \{ \Gamma \left ( \frac{bD}{\alpha^2}\right) \right \}\nonumber \\
& =\left ( \frac{bD}{\alpha^2}-\frac{1}{2}\right) \log \left (\frac{bD}{\alpha^2}\right) - \frac{bD}{\alpha^2}+\frac{1}{2}\log(2\pi)+O(\alpha^2)\nonumber \\
& =-\frac{2bD}{\alpha^2} \log(|\alpha|) + \frac{bD}{\alpha^2}\log(b)+\frac{bD}{\alpha^2}\log(D)-\frac{1}{2}\log(b)\nonumber \\
& \hskip 0.2truein -\frac{1}{2}\log(D) +\log(|\alpha|) -\frac{bD}{\alpha^2}+\frac{1}{2}\log(2\pi)+O(\alpha^2).
\label{sum_b}
\end{align}
Similarly,
\begin{align}
\log \left \{\Gamma (b_j)\right\}&=\log\left \{\Gamma\left (  \frac{b}{\alpha^2}(1+\alpha c_j)\right )\right\}\nonumber \\
&=\left \{\frac{b}{\alpha^2}(1+\alpha c_j)-\frac{1}{2}\right\} \log \left \{ \frac{b}{\alpha^2}(1+\alpha c_j)\right\}\nonumber \\
& \hskip 1.0truein -\frac{b}{\alpha^2}(1+\alpha c_j)
+\frac{1}{2}\log(2\pi)+O(\alpha^2) \nonumber \\
&=\frac{b}{\alpha^2}\log(b)-\frac{2b}{\alpha^2}\log(|\alpha|)+\frac{b c_j}{\alpha}\log(b)-\frac{2bc_j}{\alpha}\log(|\alpha|)\nonumber \\
& \hskip 1.0truein -\frac{1}{2}\log(b)+\log(|\alpha|)+\frac{b c_j}{\alpha} \nonumber \\
&\hskip 0.2truein -\frac{b c_j^2}{2}+\frac{\alpha b c_j^3}{6} +bc_j^2-\frac{\alpha b c_j^3}{2}-\frac{\alpha c_j}{2}\nonumber \\
& \hskip 1.0truein -\frac{b}{\alpha^2}(1+\alpha c_j)
+\frac{1}{2} \log(2\pi)+O(\alpha^2).
\label{sum_b_j}
\end{align}
Summing (\ref{sum_b_j}) over $j=1, \ldots , D$ and subtracting from (\ref{sum_b}) yields (\ref{lemma1}) after some further calculations.  \hfill$\square$
\vskip 0.2truein

\noindent {\bf Lemma 2}.  \textit{As $\alpha \rightarrow 0$},
\begin{align}
\sum_{j=1}^D \left \{ \frac{b}{\alpha^2}(1+\alpha c_j) \right \}\log \left (\sum_{k=1}^D x_{ik}^\alpha \right )
&=\frac{bD}{\alpha^2}\log(D)+\frac{bD}{\alpha} \hat{\kappa}_{1,i} + \frac{bD}{2} \hat{\kappa}_{2,i}\nonumber \\
& \hskip 1.0truein +\frac{\alpha b D}{6}\hat{\kappa}_{3,i} +O(\alpha^2),
\label{lemma2}
\end{align}
\textit{where} $\hat{\kappa}_{1,i}$, $\hat{\kappa}_{2,i}$ \textit{and} $\hat{\kappa}_{3,i}$ \textit{are defined in (\ref{kappa12}) and (\ref{kappa3}).}
\medskip

\vskip 0.2truein

\noindent{\bf Proof of Lemma 2}.   As $\alpha \rightarrow 0$,
 \begin{align}
\log \left(\sum_{k=1}^D x_{ik}^\alpha \right )&=\log \left(\sum_{k=1}^D \exp(\alpha y_{ik})\right )\nonumber \\
&=\log \left (  \sum_{k=1}^D 1+\alpha y_{ik}+\frac{1}{2}\alpha^2 y_{ik}^2+\frac{1}{6} \alpha^3 y_{ik}^3+O(\alpha^4)\right )\nonumber \\
&=\log \left \{ D+\alpha D \bar{y}_{i+} +\frac{1}{2}\alpha^2 \sum_{k=1}^D y_{ik}^2 +\frac{1}{6} \alpha^3 \sum_{k=1}^D y_{ik}^3 +O(\alpha^4) \right \}\nonumber \\
&=\log (D) +\log \left \{ 1+\alpha \bar{y}_{ik} +\frac{1}{2}\alpha^2 D^{-1} \sum_{k=1}^D y_{ik}^2 +\frac{1}{6}\alpha^3D^{-1} \sum_{k=1}^D y_{ik}^3 +O(\alpha^4)  \right \}\nonumber \\
& =\log (D) +\alpha \hat{\kappa}_{1,i}+\frac{\alpha^2}{2} \hat{\kappa}_{2,i}+\frac{\alpha^3}{6}\hat{\kappa}_{3,i}+O(\alpha^4),
\label{log_sum_power}
\end{align}
where $\hat{\kappa}_{i,1}$, $\hat{\kappa}_{i,2}$ and $\hat{\kappa}_{i,3}$ are defined in (\ref{kappa12}) and (\ref{kappa3}).
Since, by (\ref{sum_c}),
\[
\sum_{j=1}^D \frac{b c_j}{\alpha} \log \left ( \sum_{k=1}^D x_{ik}^\alpha \right )=0,
\]
(\ref{lemma2}) follows easily when we multiply (\ref{log_sum_power}) by $(bD)/\alpha^2$ and collect terms.  \hfill$\square$

\section{Discussion and Conclusions}

Our numerical and theoretical results show that some care is needed when applying a finite-dimensional family of transformations to a parametric model, and then fitting the resulting larger model by maximum likelihood.  The particular focus of study here is on using the $\alpha$-transformed Dirichlet distribution as the parametric model. The numerical results in Section 5 and the theoretical results Proposition 1 and Corollary 1 in Section 4 suggest that when the maximum likelihood estimator of $\alpha$ is close to $0$, there is a tendency for the other parameter estimates to be rather large.  In some of these cases with $\hat{\alpha}$ small one may be better off using the log-ratio transformation, provided there are no sample data components which are zero or very close to zero.

\subsection*{Acknowledgements}
This work was partially supported by EPSRC grant EP/K022547/1, for which we are grateful.
Partial results from this research were obtained when the second author was a PhD student at the University of Nottingham.

\end{document}